\documentclass[twoside,12pt]{article}
\usepackage[]{amsmath,amsthm,amssymb}
\usepackage[latin1]{inputenc}

\begin{document}

\title{{\Large\bf  Discrete  index  transforms with Bessel and modified Bessel functions}}

\author{Semyon  YAKUBOVICH}
\maketitle

\markboth{\rm \centerline{ Semyon   YAKUBOVICH}}{}
\markright{\rm \centerline{  DISCRETE TRANSFORMS }}

\begin{abstract} {\noindent Discrete analogues of the index  transforms, involving Bessel and the modified Bessel functions are introduced and investigated. The corresponding inversion theorems for   suitable classes  of functions and sequences   are established. }

\end{abstract}
\vspace{4mm}



\vspace{4mm}

\section {Introduction and preliminary results}

Our goal in this paper is to investigate the mapping properties and prove inversion formulas for the following six  transformations between suitable sequences $\{a_n\}_{n\ge 1}$ and functions $f$ in terms of the series and integrals,  respectively, which are associated with Bessel and  modified Bessel functions $J_\mu(z), I_\mu(z), K_\mu(z)$ (cf. [1],  Ch. 10), namely,

$$f(x)=  e^{-x/2} \sum_{n=1}^\infty a_n  {\rm Re} \left[ I_{in}\left({x\over 2}\right) \right],\quad x > 0,\eqno(1.1)$$

$$ a_n=   \int_0^\infty   {\rm Re} \left[ I_{in} \left({x\over 2}\right) \right] f(x) e^{-x/2} dx,\quad n \in \mathbb{N}_0,\eqno(1.2)$$

$$f(x)=  \sum_{n=0}^\infty {a_n \over \cosh(\pi n/2)}\   {\rm Re} \left[ J_{in}\left(2\sqrt {2x} \right) \right] K_{in}\left(2\sqrt {2x}\right),\quad x > 0,\eqno(1.3)$$

$$ a_n=  {1\over \cosh(\pi n/2)} \int_0^\infty   {\rm Re} \left[ J_{in}\left(2\sqrt {2x}\right) \right] K_{in}\left(2\sqrt {2x}\right) f(x) dx,\quad n \in \mathbb{N}_0,\eqno(1.4)$$

$$f(x)=  \sum_{n=1}^\infty {a_n  \over \sinh(\pi n/2)} \  {\rm Im} \left[ J_{in}\left(2\sqrt {2x}\right) \right]K_{in}\left(2\sqrt {2x}\right),\quad x > 0,\eqno(1.5)$$

$$ a_n=  {1\over \sinh (\pi n/2)} \int_0^\infty   {\rm Im} \left[ J_{in} \left(2\sqrt {2x}\right) \right] K_{in}\left(2\sqrt {2x}\right) f(x) dx,\quad n \in \mathbb{N}.\eqno(1.6)$$
Here  $i$ is the imaginary unit and ${\rm Re}, \ {\rm Im}$ denote the  real and imaginary parts of a complex-valued function.   We call transformations (1.1)-(1.6) the discrete index transforms  (cf. [3]). For instance, continuous analogues of transformations (1.1), (1.2) were considered  by the author in [4].   Bessel function $J_\nu(z)\ z,\nu \in \mathbb{C}$ of the first  kind is a solution of the Bessel differential equation 
$$  z^2{d^2u\over dz^2}  + z{du\over dz} + (z^2- \nu^2)u = 0.\eqno(1.7)$$
This function has  the following asymptotic behavior at infinity and near the origin
$$ J_\nu(z) = \sqrt{2\over \pi z} \cos \left( z- {\pi\over 4} (2\nu+1)\right)  [1+ O(1/z)], \ z \to \infty,\   |\arg z| <  \pi,\eqno(1.8)$$
$$J_\nu(z) = O( z^{\nu} ), \ z \to 0,\eqno(1.9)$$
The modified Bessel functions of the first and second kind $I_\nu(z),\ K_\nu(z),\ z,\nu \in \mathbb{C}$, in turn,  are solutions of the modified Bessel differential equation 
$$  z^2{d^2u\over dz^2}  + z{du\over dz} - (z^2+ \nu^2)u = 0,\eqno(1.10)$$
having the corresponding asymptotic behavior 
$$I_\nu(z) = O\left( |z|^{{\rm Re}\nu} \right), \ z \to 0,\eqno(1.11)$$
$$I_\nu(z) = O\left( {e^z\over \sqrt{2\pi z} } \right), \ z \to \infty,   \ - {\pi\over 2} < \arg z <  {3\pi\over 2},\eqno(1.12)$$
$$K_\nu(z) = O\left( |z|^{-|{\rm Re}\nu|} \right), \ z \to 0,\ \nu\neq 0,\ K_0(z) = O\left(\log(|z|)\right),  \ z \to 0,\eqno(1.13)$$
$$K_\nu(z) = O\left( \sqrt{\pi\over 2z}\  e^{-z} \right), \ z \to \infty, \  | \arg z| <  {3\pi\over 2}.\eqno(1.14)$$
These functions are related by the equality 
$$K_\nu(z)= {\pi\over 2\sin(\pi\nu)} \left[ I_{-\nu}(z)- I_\nu(z)\right].\eqno(1.15)$$
The modified Bessel function of the first kind satisfies the  inequality (cf. [4], p. 138) 

$$\left|I_{i\tau}(x)\right| \le I_0(x)  \left({\sinh(\pi \tau)\over \pi \tau} \right)^{1/2},\quad x >0,\ \tau \in \mathbb{R}.\eqno(1.16)$$
Meanwhile, the modified Bessel function of the second kind obeys the estimate [3], p. 219

$$\left| K_{i\tau}(x)\right| \le A {x^{-1/4} \over \sqrt{\sinh(\pi\tau)}},\quad x,\ \tau > 0,\eqno(1.17)$$
where $A > 0$ is an absolute constant.  The Mellin-Barnes integral representation for the kernel in (1.1), (1.2) is established in [4].  We have

$${\sqrt\pi\over \cosh(\pi \tau) } e^{-x/2} {\rm Re} \left[ I_{i\tau}\left({x\over 2}\right) \right] = {1\over 2\pi i}  \int_{\gamma-i\infty}^{\gamma +i\infty}  \frac{\Gamma\left(1/2- s \right)  \Gamma\left(s+ i\tau\right)\Gamma\left(s-i\tau\right)}{ \Gamma (s) \Gamma(1-s)} x^{-s} ds,\eqno(1.18) $$
where $x >0, \tau \in \mathbb{R}$, $\Gamma(z)$  is the Euler gamma function (cf. [1],  Ch. 5), and the contour is  the vertical straight line $s= \gamma + it,\ 0 < \gamma < 1/2,\ t \in \mathbb{R}$ in complex plane, separating the left-hand simple poles from the right-hand  ones  in the numerator of  the integrand.  In the meantime, applying the Stieltjes transform to both sides of (1.18) and changing the order of integration by Fubini's theorem via the estimate

$$ \int_0^\infty {t^{-\gamma} dt  \over x+t}  \int_{\gamma-i\infty}^{\gamma +i\infty}  \left| \frac{\Gamma\left(1/2- s \right)  \Gamma\left(s+ i\tau\right)\Gamma\left(s-i\tau\right)}{ \Gamma (s) \Gamma(1-s)}   ds \right| < \infty,$$
we find, employing Entries 8.4.2.5 and 8.4.23.5 in [2], Vol. III, the following integral representation of the modified Bessel function of the second kind (1.15)

 $$K_{i\tau} \left({x\over 2} \right) =  e^{-x/2} \int_0^\infty  {e^{-t/2}\over x+t}  {\rm Re} \left[ I_{i\tau}\left({t\over 2}\right) \right] dt.\eqno(1.19)$$
But since (cf. Entry 2.16.6.1 in [2], Vol. II)

$$\int_0^\infty e^{-x\cosh(u)} K_{in}(x) dx = {\pi \sin(nu)\over \sinh(u) \sinh(\pi n) },\quad n \in \mathbb{N},\ u \in \mathbb{R},\eqno(1.20)$$
we derive from (1.19)

$$ {\pi \sin(nu)\over \sinh(u) \sinh(\pi n) } =   \int_0^\infty e^{-x (1+\cosh(u)) }  \int_0^\infty  {e^{-t/2}\over 2x+t}  {\rm Re} \left[ I_{i n}\left({t\over 2}\right) \right] dt dx.\eqno(1.21)$$
The integral with respect to $x$ is calculated in [2], Vol. I, Entry 2.3.4.3 in terms of the upper incomplete gamma function $\Gamma(\nu, x)$ (cf. [1], Ch. 8), 

$$ \Gamma(\nu, x) = \int_x^\infty t^{\nu-1} e^{-t} dt,$$
and we have 

$$\int_0^\infty {e^{-x (1+\cosh(u)) }\over 2x +t} dx =  {1\over 2} \  e^{ t \cosh^2(u/2)} \Gamma \left(0, t \cosh^2 \left({u\over 2}\right)\right).\eqno(1.22)$$  
Thus we obtain from (1.22) the value of the integral

$$\int_0^\infty  e^{ t ( \cosh^2(u/2)- 1/2)} \Gamma \left(0, t \cosh^2 \left({u\over 2}\right)\right)  {\rm Re} \left[ I_{in}\left({t\over 2}\right) \right] dt =  {2 \pi \sin(nu)\over \sinh(u) \sinh(\pi n) }.\eqno(1.23)$$
We note that the interchange of the order of integration in (1.22) is guaranteed by Fubini's theorem, owing to  the estimate

$$\int_0^\infty e^{-x (1+\cosh(u)) }  \int_0^\infty  {e^{-t/2}\over 2x+t}  \left| I_{i n}\left({t\over 2}\right) \right| dt dx$$

$$ \le {1\over 2\sqrt 2} \int_0^\infty e^{-x (1+\cosh(u)) } {dx\over \sqrt x} \int_0^1  {e^{-t/2}\over \sqrt t}  \left| I_{i n}\left({t\over 2}\right) \right| dt $$

$$+ \int_0^\infty e^{-x (1+\cosh(u)) } dx  \int_1^\infty   {e^{-t/2}\over  t}  \left| I_{i n}\left({t\over 2}\right) \right| dt < \infty,\eqno(1.24)$$
and asymptotic formulas (1.11), (1.12) for the modified Bessel function of the first kind.   The Mellin-Barnes integrals for  kernels in (1.3)-(1.6) are given in [2], Vol. III, Entry 8.4.23.11

 $$ {1\over \cosh(\pi \tau/2)}\   {\rm Re} \left[ J_{i\tau}\left(2\sqrt {2x} \right) \right] K_{i\tau}\left(2\sqrt {2x}\right) $$
 
 $$= {1\over 16\pi \sqrt \pi i}  \int_{\gamma-i\infty}^{\gamma +i\infty}  \frac{\Gamma\left((1+s)/2 \right)  \Gamma\left((s+ i\tau)/2\right)\Gamma\left((s-i\tau)/2\right)}{ \Gamma(1-s/2)} x^{-s} ds,\eqno(1.25) $$
$$ {1\over \sinh(\pi \tau/2)}\   {\rm Im} \left[ J_{i\tau}\left(2\sqrt {2x} \right) \right] K_{i\tau}\left(2\sqrt {2x}\right) $$
 
 $$=- {1\over 16\pi \sqrt \pi i}  \int_{\gamma-i\infty}^{\gamma +i\infty}  \frac{\Gamma\left(s/2 \right)  \Gamma\left((s+ i\tau)/2\right)\Gamma\left((s-i\tau)/2\right)}{ \Gamma((1-s)/2)} x^{-s} ds,\eqno(1.26) $$
where $x, \gamma > 0,\ \tau \in \mathbb{R}$.  Hence, we employ Entries 8.4.5.1 and 8.4.23.1 in [2], Vol. III to find from (1.25) 

$$  {1\over \cosh(\pi \tau/2)}\ \int_0^\infty \sin (xt)   {\rm Re} \left[ J_{i\tau}\left(2\sqrt {2t} \right) \right] K_{i\tau}\left(2\sqrt {2t}\right) dt $$
 
 $$= {1\over 16\pi  i}  \int_{\gamma-i\infty}^{\gamma +i\infty}  \Gamma\left({s+ i\tau\over 2}\right)\Gamma\left({s-i\tau\over 2}\right)  2^{-s} x^{s-1}  ds = {1\over 2 x} K_{i\tau} \left({4\over x}\right).\eqno(1.27)$$
The interchange of the order of integration can be justified via Fubini's theorem and Stirling's asymptotic formula for the gamma function [1], Ch. 5, splitting the integral by $t$ over $(0,1)$ and $(1,\infty)$ and choosing $\gamma \in (0,1)$ and $\gamma > 1$, correspondingly.  Then with (1.20), Entry 2.16.22.9 in [2], Vol. II and asymptotic formulas (1.8), (1.9), (1.13), (1.14) we derive the value of the integral $(n \in \mathbb{N},\ u \in \mathbb{R})$

$$\int_0^\infty   {\rm Im} \left[K_0\left(4 e^{\pi i/4} \cosh^{1/2}( u)  \sqrt t\right) \right]   {\rm Re} \left[ J_{in}\left(2\sqrt {2t} \right) \right] K_{in }\left(2\sqrt {2t}\right) dt $$

$$= - { \pi \sin(nu)\over 32 \sinh(u) \sinh(\pi n/2) }.\eqno(1.28)$$
In the same manner we establish the value of the integral 

$$\int_0^\infty {\rm Re} \left[K_0\left(4 e^{\pi i/4} \cosh^{1/2}( u)  \sqrt t\right) \right]  {\rm Im} \left[ J_{in}\left(2\sqrt {2t} \right) \right] K_{in }\left(2\sqrt {2t}\right) dt $$

$$=-  {\pi \sin(nu)\over 32 \sinh(u) \cosh(\pi n/2) }.\eqno(1.29)$$
In the sequel we will provide existence conditions for discrete transformations (1.1)-(1.6) and establish their inversion formulas for suitable sequences and functions.  To do this,  we will employ, in particular,  classical Fourier series for Lipschitz functions.

\section{Inversion theorems} 

We begin with

{\bf Theorem 1}. {\it   Let a sequence $ \{a_n\}_{n\in \mathbb{N}} $ satisfy the condition  

$$\sum_{n=1}^\infty  |a_n| {e^{\pi n/2}\over \sqrt n}  < \infty.\eqno(2.1)$$
Then the discrete transformation $(1.1)$ can be inverted by the formula

$$a_n =  {1\over \pi^2}  \sinh(\pi n) \int_0^\infty   \Phi_n(x) f(x) dx,\ n \in \mathbb{N}_0,\eqno(2.2)$$
where the kernel $\Phi_n(x)$ is defined by 

$$\Phi_n(x) = e^{-x/2}  \int_0^\pi e^{ x \cosh^2(u/2)} \Gamma \left(0, x \cosh^2 \left({u\over 2}\right)\right) \sinh(u) \sin(nu) du,\ x >0,\ n \in \mathbb{N}_0,\eqno(2.3)$$
and  integral  $(2.2)$ converges absolutely. }

\begin{proof}  Indeed, substituting (1.1) and (2.3) on the right-hand side of (2.2), we appeal to the estimate, using inequality (1.16), integral (1.22)  and condition (2.1), to deduce

$$  \int_0^\infty   \left| \Phi_n(x) f(x)\right| dx \le \int_0^\infty    \int_0^\pi e^{ x ( \cosh^2(u/2)- 1)} \Gamma \left(0, x \cosh^2 \left({u\over 2}\right)\right) \sinh(u) $$

$$\times  \sum_{m=1}^\infty | a_m | \left|  {\rm Re} \left[ I_{im}\left({x\over 2}\right) \right] \right|  dx du \le 2 \int_0^\pi \sinh(u) \int_0^\infty  e^{-x/2}  I_{0}\left({x\over 2}\right) \int_0^\infty {e^{-y (1+\cosh(u)) }\over 2y +x} dy dx du $$

$$\times  \sum_{m=1}^\infty  |a_m| {e^{\pi m/2}\over \sqrt m }  \le \sum_{m=1}^\infty  |a_m| {e^{\pi m/2}\over \sqrt m }  \bigg[  {1\over \sqrt 2}  \int_0^\pi \sinh(u) \int_0^1  e^{-x/2}  I_{0}\left({x\over 2}\right) {dx\over \sqrt x} \bigg.$$

$$\bigg. \times \int_0^\infty e^{-y (1+\cosh(u)) }{ dy du\over \sqrt y}    + 2 \int_0^\pi {\sinh(u)\over 1+\cosh(u)} du \int_1^\infty  e^{-x/2}  I_{0}\left({x\over 2}\right) {dx\over x}   \bigg]  < \infty.$$
 Consequently, interchanging the order of integration and summation, we recall (1.23) to obtain 

$${1\over \pi^2}  \sinh(\pi n) \int_0^\infty   \Phi_n(x) f(x) dx = {2\over \pi}  \sinh(\pi n) \sum_{m=1}^\infty  {a_m \over \sinh(\pi m) }\int_0^\pi \sin(nu) \sin(mu) du = a_n.$$
Theorem 1 is proved. 
 
 \end{proof}
 
 The discrete transformation (1.2) can be inverted by the following theorem.

{\bf Theorem 2}.   {\it Let $f$ be a complex-valued function on $\mathbb{R}_+$ which is represented by the integral 

$$f(x) =   \int_{-\pi}^\pi  e^{ x \left(\cosh^2(u/2)- 1/2\right)} \Gamma \left(0, x \cosh^2 \left({u\over 2}\right)\right)  \psi(u)\sinh(u)  du,\quad x >0,\eqno(2.4)$$ 
where  $\psi$ is a  $2\pi$-periodic function, satisfying the Lipschitz condition on $[-\pi, \pi]$, i.e.

$$\left| \psi(u) - \psi(v)\right| \le C |u-v|, \quad  \forall \  u, v \in  [-\pi, \pi],\eqno(2.5)$$
where $C >0$ is an absolute constant.  Then for all $x >0$ the following inversion formula for  transformation $(1.2)$  holds

$$ f(x)  =   {1\over  \pi^2} \sum_{n=1}^\infty  \sinh(\pi n)  \Phi_n (x) a_n,\eqno(2.6)$$
where $\Phi_n$ is defined by $(2.3)$.}

\begin{proof}    Plugging the right-hand side of the representation (2.4) in (1.2), we interchange the order of integration,  employ  (1.23)  to obtain 

$$a_n = {2\pi\over \sinh(\pi n) }\int_{-\pi}^\pi \psi (u) \sin(nu) du.\eqno(2.7)$$
This interchange is permitted  due to estimates (1.24).   Then we substitute $a_n$ by (2.7) and $\Phi_n$ by (2.3) into the partial sum of the series (2.6) $S_N(x) $, and  it becomes 

$$ S_N(x)  = {1\over \pi} \sum_{n=1}^N   \int_{-\pi}^\pi e^{ x \left(\cosh^2(t/2)- 1/2\right)} \Gamma \left(0, x \cosh^2 \left({t\over 2}\right)\right)  \sinh(t) \sin(nt) dt $$

$$\times \int_{-\pi}^\pi   \psi(u) \sin(nu) du.\eqno(2.8)$$
Hence we   calculate  the sum in (2.8), using  the known identity 
$$ \sum_{n=0}^N  \sin(nt) \sin(nu) = {1\over 4} \left[  {\sin \left((2N+1) (u-t)/2 \right)\over \sin( (u-t) /2)}  -  {\sin \left((2N+1) (u+t)/2 \right)\over \sin( (u+t) /2)} \right],\eqno(2.9)$$
 and equality (2.8) becomes  
 
 $$ S_N(x)  =   {1\over 4 \pi} \int_{-\pi}^\pi e^{ x \left(\cosh^2(t/2)- 1/2\right)} \Gamma \left(0, x \cosh^2 \left({t\over 2}\right)\right)  \sinh(t) $$
 
 $$\times  \int_{-\pi}^\pi    \left[ \psi(u) + \psi(-u) \right] {\sin \left((2N+1) (u-t)/2 \right)\over \sin( (u-t) /2)} du dt.\eqno(2.10)$$
 Since $\psi$ is $2\pi$-periodic, we treat  the latter integral with respect to $u$ as follows 

$$  \int_{-\pi}^{\pi}   \left[ \psi(u) + \psi(-u) \right] \ {\sin \left((2N+1) (u-t)/2 \right)\over \sin( (u-t) /2)}  du $$

$$=  \int_{ t-\pi}^{t+ \pi}   \left[ \psi(u) + \psi(-u) \right] \  {\sin \left((2N+1) (u-t)/2 \right)\over \sin( (u-t) /2)}  du $$

$$=  \int_{ -\pi}^{\pi}    \left[ \psi(u+t) + \psi(-u-t) \right] \  {\sin \left((2N+1) u/2 \right)\over \sin( u /2)}  du. $$
Moreover,

$$ {1\over 2\pi} \int_{ -\pi}^{\pi}   \left[ \psi(u+t) + \psi(-u-t) \right]   \  {\sin \left((2N+1) u/2 \right)\over \sin( u /2)}  du -  \psi(t) - \psi(-t)$$

$$=  {1\over 2\pi} \int_{ -\pi}^{\pi}  \left[ \psi(u+t) + \psi(-u-t) - \psi(t) - \psi(-t) \right]  \  {\sin \left((2N+1) u/2 \right)\over \sin( u /2)}  du.$$
When  $u+t > \pi$ or  $u+t < -\pi$ then we interpret  the value  $\psi(u+t)$ by  formulas

$$\psi(u+t)- \psi(t)= \psi(u+t-2\pi)- \psi(t - 2\pi),$$ 

$$\psi(u+t)- \psi(t) = \psi(u+t+ 2\pi)- \psi(t +2\pi),$$ 
respectively.     Then   due to the Lipschitz condition (2.5) we have the uniform estimate for any $t \in [-\pi,\pi]$

$${\left|  \psi(u+t) + \psi(-u-t) - \psi(t) - \psi(-t)\right| \over | \sin( u /2) |}  \le 2 C \left| {u\over \sin( u /2)} \right|.$$
Therefore,  owing to the Riemann-Lebesgue lemma

$$\lim_{N\to \infty } {1\over 2\pi} \int_{ -\pi}^{\pi}  \left[  \psi(u+t) + \psi(-u-t) - \psi(t) - \psi(-t) \right]  \  {\sin \left((2N+1) u/2 \right)\over \sin( u /2)}  du =  0\eqno(2.11)$$
for all $ t\in [-\pi,\pi].$    Moreover, since  (see (1.22)) 
$$ \int_{-\pi}^\pi   e^{ x \left(\cosh^2(t/2)- 1/2\right)} \Gamma \left(0, x \cosh^2 \left({t\over 2}\right)\right) \left| \sinh(t) \right| $$

$$\times  \int_{ -\pi}^{\pi} \left| \left[  \psi(u+t) + \psi(-u-t) - \psi(t) - \psi(-t)   \right]\  {\sin \left((2N+1) u/2 \right)\over \sin( u /2)}  \right| du dt$$

$$ \le  {B e^{-x/2} \over \sqrt x} \int_{ -\pi}^{\pi}  {\left| \sinh(t) \right| dt \over (1+\cosh(t) )^{1/2} } \int_{ -\pi}^{\pi}   \left| {u\over \sin( u /2)} \right| du < \infty,$$
where $B  >0$ is a constant.   Therefore  via  the dominated convergence theorem it is possible to pass to the limit when $N \to \infty$ under the  integral sign, and  we derive from (2.11)

$$  \lim_{N \to \infty}    \int_{-\pi}^\pi   e^{ x \left(\cosh^2(t/2)- 1/2\right)} \Gamma \left(0, x \cosh^2 \left({t\over 2}\right)\right) \sinh(t) $$

$$ \times \int_{ -\pi}^{\pi}  \left[ \psi(u+t) + \psi(-u-t) - \psi(t) - \psi(-t)  \right]  {\sin \left((2N+1) u/2 \right)\over \sin( u /2)}  du dt  = 0.$$
Hence, combining with (2.10),  we obtain  due to  the definition of $f$

$$ \lim_{N \to \infty}  S_N(x) =  {1\over 2}  \int_{-\pi }^\pi    e^{ x \left(\cosh^2(t/2)- 1/2\right)} \Gamma \left(0, x \cosh^2 \left({t\over 2}\right)\right) \left[ \psi (t) - \psi (-t) \right] \sinh(t) dt = f(x),$$
where the integral (2.4) converges since $\psi \in C[0,\pi]$.  Thus we established  (2.6), completing the proof of Theorem 2.
 
\end{proof} 

For the discrete Re-transform (1.3) with the product of Bessel functions we have the following result.

{\bf Theorem 3}. {\it   Let a sequence $ a= \{a_n\}_{n\in \mathbb{N}} \in l_1,$ i.e.  satisfy the condition  

$$||a||_1= \sum_{n=1}^\infty  |a_n|   < \infty.\eqno(2.12)$$
Then the discrete transformation $(1.3)$ has the inversion formula 

$$a_n = - {32\over \pi^2}  \sinh(\pi n) \int_0^\infty   \Psi_n(x) f(x) dx,\ n \in \mathbb{N}_0,\eqno(2.13)$$
where the kernel $\Psi_n(x)$ is defined by 

$$\Psi_n(x) =   \int_0^\pi  {\rm Im} \left[K_0\left(4 e^{\pi i/4} \cosh^{1/2}( u)  \sqrt x\right) \right]  \sinh(u) \sin(nu)    du,\ x >0,\ n \in \mathbb{N}_0,\eqno(2.14)$$
and  integral  $(2.13)$ converges absolutely. }

\begin{proof}  In fact, doing in the same manner as in the proof of Theorem 1, we substitute (1.3) and (2.14) on the right-hand side of (2.13) and interchange the order of integration and summation.  Then we obtain, using integral (1.28), the equalities

$$- {32\over \pi^2}  \sinh(\pi n) \int_0^\infty   \Psi_n(x) f(x) dx = {2\over \pi}  \sinh(\pi n)  \sum_{m=1}^\infty  {a_m \over \sinh(\pi m)}  \int_0^\pi \sin(nu) \sin(mu) du = a_n,$$
proving (2.13).  To justify the interchange of the order of integration and summation we recall Fubini's theorem by virtue of the estimate 

$$\int_0^\infty  \left|  \Psi_n(x) f(x) \right| dx \le  \int_0^\infty    I_{0}\left(2\sqrt {2x} \right) K^2_{0}\left(2\sqrt {2x}\right) dx $$

$$\times \int_0^\pi   \sinh(u) du  \bigg[ |a_0| + \sum_{m=1}^\infty{ | a_m | \sqrt{\sinh(\pi m)} \over \sqrt {\pi m} \cosh(\pi m/2)} \bigg] $$

$$\le ||a||_1  \left(  \cosh(\pi) -1\right) \int_0^\infty    I_{0}\left(2\sqrt {2x} \right) K^2_{0}\left(2\sqrt {2x}\right) dx < \infty,\eqno(2.15)$$
which is based, in turn, on definitions of Bessel functions, inequalities (1.16), $|K_{i\tau}(x)| \le K_0(x),\ x >0, \tau \in \mathbb{R}$ and asymptotic formulas (1.11)-(1.14).  Theorem 3 is proved. 

\end{proof}

For the companion (1.4) the inversion formula in terms of the series is given by the following theorem.

{\bf Theorem 4}.   {\it For a class of functions $f$ being represented  by the integral 

$$f(x) =    - 32\int_{-\pi}^\pi  \psi(u)\sinh(u)\  {\rm Im} \left[K_0\left(4 e^{\pi i/4} \cosh^{1/2}( u)  \sqrt x\right) \right]    du,\quad x >0,\eqno(2.16)$$ 
where  $\psi$ is a  $2\pi$-periodic function, satisfying the Lipschitz condition $(2.5)$ on $[-\pi, \pi]$,  the following inversion formula   holds

$$ f(x)  =  - {32\over  \pi^2} \sum_{n=1}^\infty  \sinh(\pi n)  \Psi_n (x) a_n,\eqno(2.17)$$
where $a_n$ is given reciprocally by $(1.4)$ and $\Psi_n$ is defined by $(2.14)$.}

\begin{proof}   Doing in the same manner, we substitute the right-hand side of (2.16) in (1.4), we interchange the order of integration and  employ  (1.28)  to obtain 

$$a_n = {2 \pi\over \sinh(\pi n) }\int_{-\pi}^\pi \psi (u) \sin(nu) du.\eqno(2.18)$$
This interchange is allowed by Fubini's theorem due to the absolute convergence of the iterated integral (cf. (2.15)).   Then the partial sum of the series (2.17) has the form 

$$ S_N(x)  = - {32\over \pi} \sum_{n=1}^N  \int_{-\pi}^\pi \sinh(t) \sin(nt) {\rm Im} \left[K_0\left(4 e^{\pi i/4} \cosh^{1/2}( t)  \sqrt x\right) \right]   dt$$

$$\times  \int_{-\pi}^\pi   \psi(u) \sin(nu) du.\eqno(2.19)$$
Hence, recalling  identity (2.9), the latter equality (2.19) turns to be  
 $$ S_N(x)  =  - {8\over  \pi} \int_{-\pi}^\pi    \sinh(t)  {\rm Im} \left[K_0\left(4 e^{\pi i/4} \cosh^{1/2}( t)  \sqrt x\right) \right]  $$
 
 $$\times  \int_{-\pi}^\pi    \left[ \psi(u) + \psi(-u) \right] {\sin \left((2N+1) (u-t)/2 \right)\over \sin( (u-t) /2)} du dt.$$
Now the same scheme as in the proof of Theorem 3 drives us at the equality 

$$ \lim_{N \to \infty}  S_N(x) =  - 16  \int_{-\pi }^\pi \left[ \psi (t) - \psi (-t) \right]  \sinh(t) {\rm Im} \left[K_0\left(4 e^{\pi i/4} \cosh^{1/2}( t)  \sqrt x\right) \right] dt = f(x),$$
proving the inversion formula (2.17).

\end{proof} 

The same scheme can be applied to invert discrete $Im$-transformations (1.5), (1.6).  We will state the corresponding theorems, leaving their proofs via (1.29) to interested readers.

{\bf Theorem 5}. {\it   Let a sequence $ a= \{a_n\}_{n\in \mathbb{N}} \in l_1.$ Then the discrete transformation $(1.5)$ has the inversion formula 

$$a_n = - {32\over \pi^2}  \sinh(\pi n) \int_0^\infty   \Omega_n(x) f(x) dx,\ n \in \mathbb{N}_0,\eqno(2.20)$$
where the kernel $\Omega_n(x)$ is defined by 

$$\Omega_n(x) =    \int_0^\pi   {\rm Re} \left[K_0\left(4 e^{\pi i/4} \cosh^{1/2}( u)  \sqrt x\right) \right]  \sinh(u) \sin(nu) du,\ x >0,\ n \in \mathbb{N}_0,\eqno(2.21)$$
and  integral  $(2.20)$ converges absolutely. }

{\bf Theorem 6}.   {\it For a class of functions $f$ such that  

$$f(x) =   -32 \int_{-\pi}^\pi   {\rm Re} \left[K_0\left(4 e^{\pi i/4} \cosh^{1/2}( u)  \sqrt x\right) \right]   \psi(u)\sinh(u)  du,\quad x > 0$$
 with  $2\pi$-periodic function $\psi$, satisfying the Lipschitz condition $(2.5)$ on $[-\pi, \pi]$,  the following inversion formula   holds

$$ f(x)  =  -  {32\over  \pi^2} \sum_{n=1}^\infty  \sinh(\pi n)  \Omega_n (x) a_n,$$
where $a_n$ is given  by $(1.6)$ and $\Omega_n$ is defined by $(2.21)$.}

\bigskip
\centerline{{\bf Acknowledgments}}
\bigskip

\noindent The work was partially supported by CMUP, which is financed by national funds through FCT (Portugal)  under the project with reference UIDB/00144/2020.

\bigskip
\centerline{{\bf References}}
\bigskip
\baselineskip=12pt
\medskip
\begin{enumerate}

\item[{\bf 1.}\ ] NIST Digital Library of Mathematical Functions. http://dlmf.nist.gov/, Release 1.0.17 of 2017-12-22. F. W. J. Olver, A. B. Olde Daalhuis, D. W. Lozier, B. I. Schneider, R. F. Boisvert, C. W. Clark, B. R. Miller and B. V. Saunders, eds.

\item[{\bf 2.}\ ] A.P. Prudnikov, Yu.A. Brychkov and O.I. Marichev, {\it Integrals and Series}. Vol. I: {\it Elementary
Functions}, Vol. II: {\it Special Functions}, Gordon and Breach, New York and London, 1986, Vol. III : {\it More special functions},  Gordon and Breach, New York and London,  1990.

\item[{\bf 3.}\ ] S. Yakubovich, {\it Index Transforms}, World Scientific Publishing Company, Singapore, New Jersey, London and
Hong Kong, 1996.

\item[{\bf 4.}\ ] S. Yakubovich,  New index transforms of the Lebedev-Skalskaya type.  {\it Integral Transforms and Special Functions}, {\bf 27} (2016), N 2, 137-152.

\end{enumerate}

\vspace{5mm}

\noindent S.Yakubovich\\
Department of  Mathematics,\\
Faculty of Sciences,\\
University of Porto,\\
Campo Alegre st., 687\\
4169-007 Porto\\
Portugal\\
E-Mail: syakubov@fc.up.pt\\

\end{document}